\documentclass[reqno,a4paper]{amsart}
\usepackage[utf8]{inputenc}
\usepackage[active]{srcltx}
\usepackage{verbatim}
\usepackage{eurosans}

\allowdisplaybreaks

\usepackage{mathrsfs}
\usepackage{newsymbol}
\let\mathcal \undefined
\def\mathcal{\mathscr}
\let\emptyset \undefined
\let\ge       \undefined
\let\le       \undefined
\let\geq       \undefined
\let\leq       \undefined
\newsymbol\le          1336  \let\leq\leq
\newsymbol\ge          133E  \let\geq\geq
\newsymbol\leq          1336  \let\leq\leq
\newsymbol\geq          133E  \let\geq\geq
\newsymbol\emptyset    203F
\newsymbol\notle       230A
\newsymbol\notge       230B

\usepackage{graphicx,color,epsfig,latexsym,longtable,amsmath,amscd,amsxtra,amssymb,amsxtra,exscale,xy,mathrsfs,wrapfig,bbm}\xyoption{all}

\usepackage[fleqn,tbtags]{mathtools}
\mathtoolsset{showonlyrefs,showmanualtags}

\usepackage[plain]{fancyref}
\newcommand*{\fancyrefthmlabelprefix} {thm}
\frefformat{plain}{\fancyrefthmlabelprefix}{theorem\fancyrefdefaultspacing#1}
\Frefformat{plain}{\fancyrefthmlabelprefix}{Theorem\fancyrefdefaultspacing#1}
\newcommand*{\fancyrefcorlabelprefix} {cor}
\frefformat{plain}{\fancyrefcorlabelprefix}{corollary\fancyrefdefaultspacing#1}
\Frefformat{plain}{\fancyrefcorlabelprefix}{Corollary\fancyrefdefaultspacing#1}
\newcommand*{\fancyreflemlabelprefix} {lem} 
\frefformat{plain}{\fancyreflemlabelprefix}{lemma\fancyrefdefaultspacing#1}
\Frefformat{plain}{\fancyreflemlabelprefix}{Lemma\fancyrefdefaultspacing#1}
\newcommand*{\fancyrefproplabelprefix}{prop}
\frefformat{plain}{\fancyrefproplabelprefix}{proposition\fancyrefdefaultspacing#1}
\Frefformat{plain}{\fancyrefproplabelprefix}{Proposition\fancyrefdefaultspacing#1}
\newcommand*{\fancyrefremlabelprefix} {rem}
\frefformat{plain}{\fancyrefremlabelprefix}{remark\fancyrefdefaultspacing#1}
\Frefformat{plain}{\fancyrefremlabelprefix}{Remark\fancyrefdefaultspacing#1}
\newcommand*{\fancyrefapplabelprefix} {app}
\frefformat{plain}{\fancyrefapplabelprefix}{application\fancyrefdefaultspacing#1}
\Frefformat{plain}{\fancyrefapplabelprefix}{Application\fancyrefdefaultspacing#1}
\newcommand*{\fancyrefdeflabelprefix} {def}
\frefformat{plain}{\fancyrefdeflabelprefix}{definition\fancyrefdefaultspacing#1}
\Frefformat{plain}{\fancyrefdeflabelprefix}{Definition\fancyrefdefaultspacing#1}
\newcommand*{\fancyrefnotlabelprefix} {not}
\frefformat{plain}{\fancyrefnotlabelprefix}{notation\fancyrefdefaultspacing#1}
\Frefformat{plain}{\fancyrefnotlabelprefix}{Notation\fancyrefdefaultspacing#1}
\newcommand*{\fancyrefexlabelprefix} {ex}
\frefformat{plain}{\fancyrefexlabelprefix}{example\fancyrefdefaultspacing(#1)}
\Frefformat{plain}{\fancyrefexlabelprefix}{Example\fancyrefdefaultspacing(#1)}
\newcommand*{\fancyrefcondlabelprefix} {cond}
\frefformat{plain}{\fancyrefcondlabelprefix}{condition\fancyrefdefaultspacing(#1)}
\Frefformat{plain}{\fancyrefcondlabelprefix}{Condition\fancyrefdefaultspacing(#1)}
\usepackage[colorlinks,linkcolor={black},
citecolor={black},urlcolor={red},]{hyperref}

\theoremstyle{plain}
\newtheorem{theorem}{Theorem}[section]         
\newtheorem{corollary}[theorem]{Corollary}     
\newtheorem{lemma}[theorem]{Lemma}             
\newtheorem{proposition}[theorem]{Proposition}  
\newtheorem{conjecture}[theorem]{Conjecture}  
\theoremstyle{remark}
     
\newtheorem{example}[theorem]{Example}
\newtheorem{definition}[theorem]{Definition}

 \numberwithin{equation}{section}

\newcommand{\SCP}[2]{ {\rm(SCP)}_{(#1,#2)}}

\renewcommand {\Re} {\textit{Re}}

\newcommand {\norm}[1] {\| #1 \|}

\newcommand {\Bignorm}[1]{\Bigl\| #1 \Bigr\|}
\newcommand {\einhalb} {\frac{1}{2}}
\newcommand {\pihalbe} {\frac{\pi\,}{2}}

\newcommand {\RR} {\mathbb R}
\newcommand {\NN} {\mathbb N}

\newcommand {\ZZ} {\mathbb Z}

\newcommand {\EE} {\mathbb E}

\newcommand {\cT} {\mathscr T}
\newcommand {\cS} {\mathscr S}

  \newcommand {\la}{\lambda}
  \newcommand {\si}{\sigma}

  \newcommand {\ga}{\gamma}
  \newcommand {\om}{\omega}
  \newcommand {\Om}{\Omega}

  \renewcommand {\O}{\Omega}

  \newcommand {\al}{\alpha}
  \newcommand {\eps}{\varepsilon}

\newcommand{\suchthat}{: \;}

\newcommand{\ugr}{_{{\rm unif\hbox{-}}\ga}}
\def\N{{\mathbb N}}
\def\Z{{\mathbb Z}}

\def\R{{\mathbb R}}
\def\C{{\mathbb C}}

\newcommand{\E}{{\mathbb E}}
\renewcommand{\P}{{\mathbb P}}

\newcommand{\B}{{\mathscr B}}

\newcommand{\beq}{\begin{equation}}
\newcommand{\eeq}{\end{equation}}
\newcommand{\bal}{\begin{aligned}}
\newcommand{\eal}{\end{aligned}}
\newcommand{\ben}{\begin{enumerate}}
\newcommand{\een}{\end{enumerate}}
\newcommand{\bit}{\begin{itemize}}
\newcommand{\eit}{\end{itemize}}

\newcommand{\bth}{\begin{theorem}}
\renewcommand{\eth}{\end{theorem}}
\newcommand{\bpr}{\begin{proposition}}
\newcommand{\epr}{\end{proposition}}
\newcommand{\ble}{\begin{lemma}}
\newcommand{\ele}{\end{lemma}}
\newcommand{\bpf}{\begin{proof}}
\newcommand{\epf}{\end{proof}}
\newcommand{\bex}{\begin{example}}
\newcommand{\eex}{\end{example}}
\newcommand{\bre}{\begin{example}}
\newcommand{\ere}{\end{example}}

\renewcommand{\Re}{\hbox{\rm Re}\,}

\newcommand{\calL}{{\mathscr L}}
\newcommand{\n}{\Vert}

\newcommand{\s}{^*}
\newcommand{\lb}{\langle}
\newcommand{\rb}{\rangle}

\newcommand{\limn}{\lim_{n\to\infty}}

\newcommand{\limj}{\lim_{j\to\infty}}

\newcommand{\sumk}{\sum_{k\geq 1}}

\begin{document}

\title[Uniformly $\ga$--radonifying families]
{Uniformly $\ga$--radonifying families of operators 
and the stochastic Weiss conjecture}

\author{Bernhard H. Haak}

\address{Institut de Math\'ematiques de Bordeaux\\Universit\'e Bordeaux 1\\351 cours de la Lib\'eration\\33405 Talence CEDEX\\France}
\email{Bernhard.Haak@math.u-bordeaux1.fr}

\author{Jan van Neerven}
\address{Delft Institute of Applied Mathematics\\
Technical University of Delft\\P.O. Box 5031\\2600 GA Delft\\The Netherlands}
\email{J.M.A.M.vanNeerven@tudelft.nl}

\thanks{The second-named author is supported by
VICI subsidy 639.033.604 of the Netherlands Organisation for Scientific Research (NWO)}


\date\today

\begin{abstract}
We introduce the notion of uniform $\gamma$--radonification of a family of
operators,
which unifies the notions of $R$--boundedness of a family of operators 
and $\ga$--radonification of an individual operator.
We study the properties of uniformly $\gamma$--radonifying families of
operators in detail and apply
our results to 
the stochastic abstract Cauchy problem 
$$dU(t) = AU(t)\,dt + B\,dW(t), \quad U(0)=0.$$
Here, $A$ is the generator 
of a strongly continuous semigroup of operators on a Banach space $E$, 
$B$ is a bounded linear operator from a separable Hilbert space $H$ into $E$, 
and $W_H$ is an $H$--cylindrical Brownian motion. 
When $A$ and $B$ are simultaneously diagonalisable, 
we prove 
that an invariant measure exists if and only if 
the family 
$$ \{\sqrt{\la}R(\la,A)B: \ \la\in S_\vartheta\}$$
is uniformly $\ga$--radonifying for some/all $0< \vartheta<\pihalbe$, where
$S_\vartheta$ is the open sector of angle $\vartheta$ in the complex plane.
This result can be viewed as a partial solution of a
stochastic version of the Weiss conjecture in linear systems theory.
\end{abstract}

\maketitle

\section{ Introduction}
In recent years it has become apparent that many classical results in operator
theory and harmonic analysis 
can be generalised from their traditional Hilbert space setting to 
Banach spaces, provided the notion of uniform boundedness is replaced
with {\em $R$--boundedness}. 
This notion appeared implicitly in the work of Bourgain \cite{Bourgain}
and was formalised by Berkson and Gillespie \cite{BG} and Cl\'ement, de Pagter,
Sukochev, and Witvliet \cite{ClementDePagterSukochevWitvliet}. 
It has accomplished remarkable progress in the theory
of parabolic evolution equations. A highlight is the
recent solution of the $L^p$-maximal
regularity problem by Weis \cite{Weis:FM}, who proved an extension of the Mihlin multiplier theorem for operator-valued multipliers taking $R$--bounded values in a UMD space $E$ and used it to deduce that the generator $A$ of a bounded analytic semigroup on a UMD space $E$ has $L^p$-maximal regularity if and only if $\la \mapsto \la(\la-A)^{-1}$ is $R$--bounded on $\C_+$. Soon, an alternative approach to the $L^p$-maximal regularity problem via $H^\infty$--calculus appeared. In the Hilbert space setting this calculus
was introduced by McIntosh \cite{McIntosh:H-infty-calc}, who characterised 
it by means of square function estimates. 
This characterisation extends to Banach spaces, provided the square 
functions are replaced with $\ga$--radonifying norms
\cite{CowlingMcIntoshDoustYagi,KaltonWeis:euclidian-structures,KaltonWeis:square-function-est,LeMerdy:square-functions}.

These developments have been documented in detail in 
the memoir by Denk, Hieber, and Pr\"uss \cite{DHP} and the lecture notes by Kunstmann and Weis \cite{KunstmannWeis:Levico}, 
where extensive references can be found.

In a parallel development,
$\ga$--radonifying norms have been used recently to extend 
the theory of stochastic integration to the Banach space
setting, first for operator-valued functions taking values in arbitrary Banach spaces \cite{vanNeervenWeis:stoch-int, vanNeervenWeis:weak-limits} and subsequently for operator-valued processes taking values in UMD spaces \cite{NVW}. 
In both papers, the role of the It\^o isometry is taken over by an isometry in terms of $\ga$--radonifying norms. Applications to stochastic evolution equations in Banach spaces have been worked out, for linear equations \cite{DNW, HNV, vanNeervenWeis:asymptotic-scp} and nonlinear equations \cite{NVW2,VeraarZimmerschied}.

Further applications of $R$--boundedness and $\ga$--radonifying norms have been given in various areas on analysis, such as harmonic analysis
\cite{ArendtBu, BG, GirardiWeis, Hytonen, HytonenWeis, KaiserWeis:wavelet-trafo},
Banach space theory \cite{DiestelJarchowTonge,
KaltonWeis:euclidian-structures, Pisier:convexbodies, TJ}, 
interpolation theory 
\cite{KaltonKunstmannWeis}, control theory \cite{Haak:Diss, HaakKunstmann:l-control-theory, LeMerdy:weiss-conj, LeMerdy:square-functions}, and noncommutative analysis \cite{JLMX, PSW}; this list of references is far from complete. 

In this paper we unify the notions of $R$--boundedness (or rather, its Gaussian analogue $\ga$-boundedness) and 
$\ga$--radonification by introducing the concept of 
uniformly $\ga$--radonifying families of operators. 
As we shall demonstrate in Sections~\ref{sec:ugr} and 
\ref{sec:compact} this is a happy marriage: 
uniformly $\ga$--radonifying families enjoy many of the good properties
both of $R$--bounded families and of $\ga$--radonifying operators.

In Section~\ref{sec:laplace} we apply our abstract results to study some 
properties of 
operator-valued Laplace transforms. 
It turns out that the Laplace 
transforms of $\ga$--radon\-ifying operators $\Phi:L^2(\R_+;H)\to E$ are 
uniformly $\ga$--radonifying both in half-planes and in sectors 
properly contained in $\C_+$. 

Natural examples of uniformly $\ga$--radonifying families of operators 
arise in the theory of stochastic evolution equations. These 
will be presented in the final Section~\ref{sec:weiss}
of the paper, where we apply our results on Laplace transforms to stochastic evolution equations. 
We prove that a necessary condition
for the existence of invariant measures for
the linear stochastic Cauchy problem 
$$ dU(t) = AU(t)\,dt + B\,dW_H(t),\quad U(0)=0,$$ 
where $W_H$ is an $H$--cylindrical Brownian motion and $B:H\to E$ is a bounded operator,
is that the family 
$$ \{\sqrt{\la}R(\la,A)B: \, \la\in S_\vartheta\}$$
should be uniformly $\ga$--radonifying for all $0< \vartheta <\pihalbe$, where
$S_\vartheta$ is the open sector of angle $\vartheta$ in the complex plane.
For simultaneously diagonalisable operators $A$ and $B$ 
we show that this condition is also sufficient. 
This result is a partial solution of a stochastic version of the 
Weiss conjecture in linear systems theory (see \cite{Weiss:conjectures} and the subsequent
work \cite{JP, LeMerdy:weiss-conj, Weiss:Admissibility-of-unbounded, Weiss:Carleson, Weiss:admiss-observation})
in which $L^2$-admissibility of the control operator
is replaced with the existence of an invariant measure. 

\section{Uniformly $\ga$--radonifying families}\label{sec:ugr}

In this section we introduce the notion of a uniformly $\ga$--radonifying family of operator and study its properties. This notion unifies the 
concepts of $R$--boundedness (or rather, $\ga$--boundedness) and 
$\ga$--radonification, which we shall discuss first.

Let $E$ and $F$ be Banach spaces. A subset $\cS$ of $\B(E,F)$ is
called {\em $R$--bounded} if there exists a constant $C\geq 0$ such
that  for all $n\geq 1$, all $x_1,\dots,x_n\in E$, and all
$S_1,\dots,S_n\in\cS$ we have  $$ \E\Big\n \sum_{k=1}^n r_k S_k
x_k\Big\n^2 \leq C^2  \E\Big\n \sum_{k=1}^n r_k  x_k\Big\n^2.$$
Here, $(r_k)_{k\geq 1}$ is a {\em Rademacher sequence}, i.e. a sequence of independent $\{-1,+1\}$-valued random variables on some probability space 
$(\O,\P)$ with the property that $\P(r_k=\pm 1) = \frac12$.  
The least admissible constant $C$ is called the {\em $R$--bound} of $\cS$,
notation $R(\cS)$. A similar definition may be given in terms of Gaussian sums, 
which leads to the concept of a {\em $\ga$--bounded} family with {\em $\ga$--bound} $\ga(\cS)$.
By a standard randomisation argument, every $R$--bounded family $\cS$
is $\ga$--bounded with $\ga(\cS)\leq R(\cS)$. If $E$ and $F$ are
Hilbert spaces, the notions of $\ga$--boundedness and $R$--boundedness
coincide with that of uniform boundedness and we have
$\ga(\cS)= R(\cS)= \sup_{S\in\cS}\n S\n$.

Throughout this paper, unless otherwise stated 
 $H$ is a separable infinite-dimens\-ional Hilbert space and $E$ is a Banach space. 
Let $(\ga_k)_{k\geq 1}$ be a {\em Gaussian sequence}, i.e., a sequence of 
independent real-valued standard Gaussian random variables on some probability 
space $(\Om, \P)$.
A linear operator $T: H\to E$
is called $\ga$--{\em radonifying} if for some orthonormal basis 
$(h_k)_{k\geq 1}$ of $H$ the sum
$$ \sumk \ga_k Th_k $$
converges in $L^2(\O;E)$. If this is the case, the sum $ \sumk \ga_k Th_k $ converges in
$E$
almost surely and in $L^p(\O;E)$ for all $1\leq p<\infty$, for every
orthonormal basis $(h_k)_{k\geq 1}$ of $H$.
The linear space of all $\ga$--radonifying operators from $H$ to $E$
is denoted by 
$\ga(H,E)$. 
Endowed with the norm
$$ \n T\n_{\ga(H,E)}^2 := \E \Big\n \sumk \ga_k Th_k\Big\n^2,$$
which is
independent of the basis $(h_k)_{k\geq 1}$,
the space $\ga(H,E)$ is a Banach space.
Furthermore, it is a two-sided operator ideal in $\B(H,E)$, the 
space of all bounded linear operators from $H$ to $E$.
For proofs and more information we refer to the review paper \cite{NeeCanberra}.

\begin{definition}\label{def:u-g-r}
A subset $\cT$ of $\B(H, E)$ is {\em uniformly $\ga$--radonifying} if 
for all orthonormal bases $(h_k)_{k\geq 1}$ of $H$ and sequences
$(T_k)_{k\geq 1}$ in $\cT$ the Gaussian sum $\sumk \ga_k T_k h_k$
converges in $L^2(\O;E)$. 
\end{definition}

It is important to note that this definition refers to {\em all} orthonormal 
bases of $H$. Evidently, this definition trivializes for finite-dimensional Hilbert spaces; it is mainly for this reason that we restrict our attention to infinite-dimensional $H$. 

By considering the constant sequence $T_k=T$ we see that every
operator $T$ in a uniformly $\ga$--radonifying subset $\cT$ of $\B(H,E)$
is $\ga$--radonifying, i.e., $\cT\subseteq\ga(H,E)$.

We begin our investigations with proving some simple permanence properties
of uniformly $\ga$--radonifying families of operators, resembling those of $R$--bounded and $\ga$--bounded families of operators. In what follows,
$\cT$ denotes a subset of $\B(H,E)$.

\begin{proposition}[Strong closure]\label{prop:closure} 
If $\cT$ is uniformly $\ga$--radonifying, then the closure
$\overline{\cT}$ in the strong operator topology of $\B(H,E)$ is
uniformly $\ga$--radonifying. 
\end{proposition}
\begin{proof}
Let $(\overline{T}_k)_{k\geq 1}$ be a sequence in
$\overline{\cT}$. Given an $\eps>0$ and an 
orthonormal basis $(h_k)_{k\geq 1}$ of $H$, choose a sequence
$(T_k)_{k\geq 1}$ in $\cT$ such that 
$\n \overline{T}_k h_k - T_k h_k \n < 2^{-k}\eps$ for all $k\geq 1$.
Then, for all $1\leq M\leq N$,
$$
\bal
\Big(\E\Big\n \sum_{k=M}^N \ga_k \overline T_k h_k\Big\n^2\Big)^\frac12
& \leq  
\Big(\E\Big\n \sum_{k=M}^N \ga_k T_k h_k\Big\n^2\Big)^\frac12
+ \Big(\E\Big\n \sum_{k=M}^N \ga_k (\overline T_k h_k -T_k h_k)\Big\n^2\Big)^\frac12
\\ & \leq \Big(\E\Big\n \sum_{k=M}^N \ga_k T_k h_k\Big\n^2\Big)^\frac12 + \sum_{k=M}^N \n \overline T_k h_k -T_k h_k \n
 \\ & \leq \Big(\E\Big\n \sum_{k=M}^N \ga_k T_k h_k\Big\n^2\Big)^\frac12 +\eps.
\eal $$
The result follows by letting $M,N\to\infty$.
\end{proof}

\begin{lemma}\label{lem:tail-estimate}
If $\cT$ is uniformly $\ga$--radonifying, then for all 
orthonormal bases $(h_k)_{k\geq 1}$ of $H$ we have
$$
\lim_{n\to\infty} \Big(\sup_{T} \ \E\Big\n \sum_{k=n}^\infty \ga_k
T_k h_k\Big\n^2\Big) =0,
$$
where the supremum is taken over all sequences $T = (T_k)_{k\geq 1}$ in $\cT$.
\end{lemma}
\begin{proof}
If the lemma was false, we could find an orthonormal basis
$(h_k)_{k\geq 1}$ of $H$, a number $\delta>0$, an increasing sequence
of indices $1\leq n_1 < N_1 < n_2 < N_2 < \dots$, and for each $j=1,2,\dots$ 
a finite set of operators $T_{n_j}, \dots, T_{N_j}\in \cT$ such that 
$$
   \E\Big\n \sum_{k=n_j}^{N_j} \ga_k T_k h_k\Big\n^2 \geq \delta^2, \quad j=1,2,\dots
$$ 
Putting $T_k := 0$ if $N_j < k <n_{j+1}$ for some $j\geq 0$ (with the
convention that $N_0=0$) we obtain a sequence $(T_k)_{k\geq 1}$ for
which the sum $\sumk \ga_k T_k h_k$ fails to converge in $L^2(\O;E)$,
and we have arrived at a contradiction.
\end{proof}

\begin{proposition}[Absolute convex hull]\label{prop:convex} 
If $\cT$ is uniformly $\ga$--radonifying, then the absolute convex hull
of $\cT$ is uniformly $\ga$--radonifying. 
\end{proposition}
\begin{proof} 
Considering real and complex parts separately and possibly replacing 
some of the $\ga_k$ by $-\ga_k$, it suffices to prove the statement in
the lemma for the convex hull of $\cT$. Furthermore, by the
contraction principle for Banach space-valued Gaussian sums,  
$\cT\cup \{0\}$ is  uniformly
$\ga$--radonifying and therefore  we may assume that $0\in\cT$.

Fix an orthonormal basis $(h_k)_{k\ge 1}$ of $H$ and an $\eps>0$, and choose $n_0\geq 1$ so large that
$$  \sup_{T} \ \E\Big\n \sum_{k=n_0}^\infty \ga_k T_k h_k\Big\n^2 < \eps^2.$$
Let $(S_k)_{k\geq 1}$ be a sequence in {\rm conv}$(\cT)$ 
and fix indices $M,N$ satisfying $n_0\leq M\leq N$. Noting that 
$$
{\rm conv}(\cT)\times \dots \times{\rm conv}(\cT) = {\rm
  conv}(\cT\times\dots\times \cT)
$$
we can find $\la_1,\dots,\la_N \in [0,1]$ with $\sum_{j=1}^N \la_j=1$
such that $S_k = \sum_{j=1}^N \la_j T_{jk}$ with $T_{jk}\in \cT$ for all
$k=M,\dots, N$.  Then,
$$
\Big(\E\Big\n \sum_{k=M}^N \ga_k S_k h_k\Big\n^2\Big)^\frac12
\leq \sum_{j=1}^N \la_j\Big(\E\Big\n  \sum_{k=M}^N \ga_k T_{jk}
    h_k\Big\n^2\Big)^\frac12 < \sum_{j=1}^N \la_j\eps =\eps.
$$
\end{proof}

Combining Propositions~\ref{prop:closure} and~\ref{prop:convex}
we obtain that the strongly closed absolutely convex hull of every
uniformly $\ga$--radonifying set is uniformly $\ga$--radonifying.
As in the case of $R$--boundedness, cf. 
\cite{DHP,  KunstmannWeis:Levico, Weis:FM}, 
this may be used to show that uniform $\ga$--radonification 
is preserved by taking integral means. In this way a number of well-known 
$R$--boundedness results can be carried over to uniformly
$\ga$--radonifying families. To give a few examples we formulate
analogues of \cite[Corollary 2.14]{KunstmannWeis:Levico}
and \cite[Propositions 2.6 and 2.8]{Weis:FM}.

\begin{proposition}\label{prop:L1}
Let $(S,\mu)$ be a $\sigma$-finite measure space and let $\cT$ be a
uniformly $\ga$--radonifying subset of $\B(H,E)$. If $f:S\to \B(H,E)$
is strongly $\mu$-measurable (in the sense that $s\mapsto f(s)h$ is strongly $\mu$-measurable for all $h\in H$) with $f(s)\in \cT$
for $\mu$-almost all $s\in S$, then for all $\phi\in L^1(S,\mu)$ the
operator $$ T_\phi h := \int_S \phi(s) f(s)h\,d\mu(s), \quad h\in H,$$
belongs to $\B(H,E)$ and the family $\{T_\phi\suchthat \ \n \phi\n_1 \leq 1\}$ 
is uniformly $\ga$--radonifying.
\end{proposition}

\begin{proposition}\label{prop:analytic}
Let $G\subseteq \C$ be an open domain and let 
$f:G\to \B(H,E)$ be an analytic function with $f(z)\in\ga(H,E)$ for
all $z\in G$.  Then for every compact subset $K\subseteq G$ the family
$\{ f(z)\suchthat z\in K\}$ is uniformly $\ga$--radonifying.
\end{proposition}

\begin{proposition}\label{prop:jordan}
Let $G\subseteq \C$ be a simply connected Jordan domain such 
that $\C\setminus G$ has nonempty interior. Let 
$f:\overline{G}\to \B(H,E)$ be uniformly bounded and strongly
continuous, analytic on $G$, and assume that $\{f(z)\suchthat z\in
\partial G\}$ is  uniformly $\ga$--radonifying. Then $\{f(z)\suchthat
\ z\in \overline{G}\}$ is  uniformly $\ga$--radonifying. 
\end{proposition}

In the situation of Proposition~\ref{prop:jordan} it follows that
$f(z)\in \ga(H,E)$ for all $z\in\overline{G}$. It will follow from
Theorem~\ref{thm:unif-gamma-compact} below that $f:\overline{G}\to
\ga(H,E)$ is continuous. 

The next lemma shows that uniform $\ga$--radonification is preserved 
under left and right multiplication.

\begin{proposition}[Ideal property]\label{prop:ideal}
Let $\tilde H$ be a separable infinite-dimensional Hil\-bert space and $\tilde E$ a Banach space.
If $\cT$ is a uniformly $\ga$--radonifying subset of $\B(H,E)$ and 
$R:\tilde H\to H$ and $S: E\to \tilde E$ are bounded operators, then
$S\cT R$ is a uniformly $\ga$--radonifying subset of $\B(\tilde H,\tilde E)$.
\end{proposition}
\begin{proof}
The left ideal property is trivial. To prove the right ideal property
we first consider the case of complex scalars. For $\tilde H=H$ the
right ideal property then follows from the well-known fact that 
the convex hull of the unitary operators on $H$  space 
are uniformly dense in the closed unit ball of $\B(H)$ (this is the 
so-called Russo-Dye theorem \cite{RusDye}; at the cost of picking up
a constant $2$ we could alternatively use the elementary fact 
that every operator in $\B(H)$
of norm less than $\frac12$ is a convex combination of at most four unitaries).

In the case of different Hilbert
spaces $H$ and $\tilde H$ 
write $TR = TRU^*\circ U$, where $U$ is an isometry from $\tilde
H$ onto $H$ and note that $\cT RU\s$ is uniformly $\ga$--radonifying
by the preceding observation.

In the case of real scalars, let $(\tilde h_k)_{k\geq 1}$ be an
orthonormal basis of $\tilde H$ and write $\sumk \ga_k T_k S\tilde h_k
= \sumk \ga_k T_k^\C S^\C \tilde h_k^\C$, where $T_k^\C$ and $S^\C$
are the complexifications of $T_k$ and $S$, and $\tilde h_k^\C = \tilde h_k+i0$.
Since $(\tilde h_k^\C)_{k\geq 1}$ is an orthonormal basis for $\tilde
H^\C$, the right-hand side converges in $L^2(\Om;E^\C)$. 
\end{proof}

We continue with a preliminary boundedness result for uniformly
$\ga$--radonifying families. It will be strengthened in Theorem
\ref{thm:bdd} below.

\begin{proposition}\label{prop:bdd}
Let $\cT$ be a uniformly $\ga$--radonifying subset of $\B(H,E)$.
Then $\cT$ is a bounded subset of $\ga(H,E)$.
\end{proposition}

\begin{proof}  
The fact that $\cT$ is contained in $\ga(H,E)$ has already been noted.
It suffices to prove that $\sup_{k\geq 1} \n T_k\n_{\ga(H,E)}<\infty$
for every sequence $(T_k)_{k\geq 1}$ in $\cT$. 

Fix an orthonormal basis $(h_k)_{k\geq 1}$ of $H$. 
By Proposition~\ref{prop:ideal} and a closed graph argument,
there exists a constant $C\geq 0$ such that for all $S\in \B(H)$ we have
$$ \E \Big\n\sumk\ga_k T_kSh_k\Big\n^2 \leq C^2\n S\n^2.$$
In particular, for all $x\s\in E\s$ this implies
$$
\bal 
\sumk \big|(Sh_k, T_k\s x\s)_H\big|^2 
& = \E \Big|\sumk \ga_k \lb T_kSh_k, x\s\rb\Big|^2
\\ & = \E \Big|\Big\lb\sumk \ga_k T_kSh_k, x\s\Big\rb\Big|^2 \leq C^2\n
       S\n^2 \n x\s\n^2.
\eal
$$
Taking $S: = h_n\otimes T_n\s x\s$ with $n\geq 1$ fixed
it follows that
$$ \n T_n\s x\s\n_H^4 \leq C^2\n T_n\s x\s\n^2 \n x\s\n^2.$$
Therefore, $\sup_{n\geq 1}\n T_n\n\leq C$. 
Next, by Lemma~\ref{lem:tail-estimate} we can find $N\geq 1$ so large that
$$ \sup_{k\geq 1} \ \E\Big\n \sum_{n=N+1}^\infty \ga_n T_k h_n\Big\n^2\leq 1.$$
But then,
$$ \sup_{k\geq 1} \n T_k\n_{\ga(H,E)} \leq 1 + CN.$$
\end{proof}

The next result explains our terminology `uniformly $\ga$--radonifying':

\begin{theorem}\label{thm:bdd}
Let $\cT$ be a uniformly $\ga$--radonifying subset of $\B(H,E)$.
Then there exists a constant $C\geq 0$ such that 
for all orthonormal bases $h=(h_k)_{k\geq 1}$ of $H$ 
and all sequences $(T_k)_{k\geq 1}$ in $\cT$ we have
$$\E\Big\n\sum_{k\geq 1} \ga_k T_k h_k\Big\n^2 \leq C^2.$$
\end{theorem}
\begin{proof}
Let  
$$
W: = \overline{{\rm abs\,conv}}(\cT)\times\overline{{\rm
    abs\,conv}}(\cT)\times\dots
$$
where ${\rm abs\,conv}(\cT)$ denotes the absolute convex hull of $\cT$
and the closure is taken in the norm of $\ga(H,E)$. Note that $W$ is absolutely convex. 
Let  $l(\cT)$ denote the vector space of all sequences $T = (T_k)_{k\geq
1}$ in $\ga(H,E)$ such that $cT\in W$ for some $c>0$. 
In view of Propositions~\ref{prop:closure},~\ref{prop:convex}, and
Proposition~\ref{prop:bdd} we may endow $l(\cT)$ with the norm
$$\n T\n_{l(\cT)} :=\inf\Big\{\frac1c\suchthat  c>0, \ cT\in W\Big\}
+  \sup_{k\geq 1} \n T_k\n_{\ga(H,E)}. $$ 
It is routine to  check that the normed space $l(\cT)$ is a Banach space. 

We fix an orthonormal basis $(\overline h_k)_{k\geq 1}$ in $H$ and 
consider the bilinear operator $\beta:  \B(H) \times l(\cT)\to L^2(\O;E)$ 
defined by $$ \beta(S,T) := \sum_{k\geq 1} \ga_k T_k S\overline h_k.$$
Note that this sum converges in $L^2(\O;E)$ thanks to Propositions
\ref{prop:closure},~\ref{prop:convex}, and \ref{prop:ideal}.
By the closed graph theorem  there is a constant $c\geq 0$ such that
for all $S\in \B(H)$ and $T\in l(\cT)$,
$$
 \E\Big\n\sum_{k\geq 1} \ga_k T_k S \overline h_k\Big\n^2 \leq c^2 
\n S\n^2 \n T\n_{l(\cT)}^2.
$$
For sequences $(T_k)_{k\geq 1}$ in
$\cT$
we have $\inf\big\{\frac1c\suchthat  c>0, \ cT\in W\big\} \leq 1$ and
consequently $\n T\n_{l(\cT)} \leq 1 + \sup_{k\geq 1} \n
T_k\n_{\ga(H,E)}\leq 1+M$, where $M := \sup_{T\in \cT} \n T\n_{\ga(H,E)}$. Hence,
$$\E\Big\n\sum_{k\geq 1} \ga_k T_k S \overline h_k\Big\n^2 \leq c^2
(1+M)^2\n S\n^2.
$$
Finally, if $(h_k)_{k\geq 1}$ is an arbitrary orthonormal basis of $H$ we let 
$U$ be the unitary operator defined by $U\overline h_k = h_k$ and
obtain, for all  sequences $(T_k)_{k\geq 1}$ in $\cT$,
$$ \E\Big\n\sum_{k\geq 1} \ga_k T_k h_k\Big\n^2
=\E\Big\n\sum_{k\geq 1} \ga_k T_k U\overline h_k\Big\n^2 \leq c^2 (1+M)^2\n U\n^2 
= c^2(1+M)^2.$$
\end{proof}

For a uniformly $\ga$--radonifying family $\cT$ in $\ga(H,E)$ we define
$$ \n \cT\n\ugr^2 := \sup_{h}\sup_{T}\ \E\Big\n\sum_{k\geq 1} \ga_k T_k h_k\Big\n^2,$$
where the fist supremum is taken over all orthonormal bases of $H$ and the
second over all sequences in $\cT$. Inspection of the proofs of
Propositions~\ref{prop:closure},~\ref{prop:convex}, and
\ref{prop:ideal} shows that we have  $$  \n \overline{\cT}\n\ugr = \n \cT\n\ugr,$$
where $\overline{\cT}$ is the strong closure of $\cT$, 
\beq\label{eq:absconv} \bal 
\n {\rm abs\,conv}({\cT})\n\ugr  & = \n \cT\n\ugr \ \ \ \hbox{(real scalars)},\\
\n {\rm abs\,conv}({\cT})\n\ugr  &\leq  2\n \cT\n\ugr \ \ \hbox{(complex scalars)},
\eal
\eeq
and
\beq\label{eq:4}
 \n {R\cT S}\n\ugr \leq \n R\n\, \n \cT\n\ugr\,\n S\n.
\eeq

Using \eqref{eq:absconv} we obtain analogous bounds for the sets discussed
in the Propositions~\ref{prop:L1},~\ref{prop:analytic}, and~\ref{prop:jordan}.

We proceed with some applications of Theorem~\ref{thm:bdd}. 
The first two results clarify the relation between uniform
$\ga$--radonification and $\ga$--boundedness.

\begin{corollary}\label{cor:gbdd}
If $\cT$ is uniformly $\ga$--radonifying, then:
\ben 
\item $\cT$ is $R$--bounded with $R(\cT) \leq \sqrt{\tfrac12\pi}\n \cT\n\ugr.$
\item $\cT$ is $\ga$--bounded with $\ga(\cT) \leq \n \cT\n\ugr;$
\een
\end{corollary}
\begin{proof}
We shall prove part (a). Since every $ R$--bounded set is $\ga$--bounded
with the same boundedness constant, the $\ga$-boundedness assertion in 
(b) follows directly from (a), but this argument produces an additional  
constant
$\sqrt{\tfrac12\pi}$. The sharper constant $1$ is obtained by noting
that for the proof of (b), Rademacher variables can be replaced by Gaussians and the first inequality in $\eqref{(1)}$ can be omitted and we may replace the role of Rademachers by Gaussians in the last step of the argument.

Fix $T_1,\dots,T_n\in \cT$ and vectors $g_1,\dots,g_n\in H$.
Let $(h_k)_{k\geq 1}$ be an orthonormal basis of $H$ 
and define $S\in \B(H)$ by $Sh_k = g_k$ for $k=1,\dots,n$ and 
$Sh_k=0$ for $k\geq n+1$. If $(r_k)_{k\geq 1}$ is a Rademacher sequence,
then
$$
\bal
\n Sh\n  & = \Big\n \sum_{k=1}^n (h,h_k)_H g_k \Big\n_H
 \leq  \sum_{k=1}^n \big|(h,h_k)_H\big|\, \n g_k\n_H 
\\ & \leq \n h\n_H \Big(\sum_{k=1}^n \n g_k\n_H^2\Big)^\frac12
=\n h\n_H \Big(\E\Big\n \sum_{k=1}^n r_k g_k \Big\n^2\Big)^\frac12.
\eal$$
Hence,  estimating Rademachers with Gaussians
and using \eqref{eq:4},
\beq\label{(1)}
\bal 
\E\Big\n \sum_{k=1}^n r_k T_k g_k \Big\n^2 
\ & \leq \tfrac12\pi\, \E\Big\n \sum_{k=1}^n \ga_k T_k S h_k \Big\n^2
\\ & 
\leq \tfrac12 \pi\n \cT\n\ugr^2 \n S\n^2
\leq \tfrac12\pi\n \cT\n\ugr^2\E\Big\n \sum_{k=1}^n r_k g_k \Big\n^2.
\eal
\eeq
\end{proof}

The next example shows that even for Hilbert spaces $E$, a
$\ga$--bounded family of operators in $\ga(H,E)$ need not be uniformly
$\ga$--radonifying.

\begin{example}
Let $(h_k)_{k\geq 1}$ be an orthonormal basis for an
infinite-dimensional Hilbert space $H$ and let $P_n$ be the orthogonal
projection onto the span of the vector $h_n$. The family
$\{P_n\suchthat  n\geq 1\}$ is uniformly bounded, hence
$\ga$--bounded, in $\B(H)$ and fails to be uniformly
$\ga$--radonifying, as is immediate by considering the sum $\sumk
\ga_k P_k h_k$.
\end{example}

The next corollary identifies $\ga$--bounded sets as the class of
`multipliers' for uniformly $\ga$--radonifying sets:

\begin{corollary} \label{cor:gbdd2}
For a subset $\cS$ of $\B(E,F)$ the following assertions are equivalent:
\ben
\item $\cS$ is $\ga$--bounded;
\item $\cS T$ is a uniformly $\ga$--radonifying subset of $\B(H,F)$ 
for every $T\in \ga(H,E)$;
\item $\cS \cT$ is a uniformly $\ga$--radonifying subset of $\B(H,F)$
for every uniformly $\ga$--radonifying subset $\cT$ of $\B(H,E)$.
\een
In the situation of {\rm (c)} we have
$ \n \cS \cT\n\ugr \leq \ga(\cS) \n\cT\n\ugr.$
\end{corollary}

\begin{proof} 
The implication (a)$\Rightarrow$(c) and the estimate 
are immediate consequences of the definitions, and the implication
(c)$\Rightarrow$(b) is trivial. To prove (b)$\Rightarrow$(a) we fix 
an orthonormal basis $(h_k)_{k\geq 1}$ in $H$ and denote by $\cS^\infty
= \cS\times\cS\times\dots$ the set of all sequences in $\cS$. By
Theorem~\ref{thm:bdd}, for each $T\in\ga(H,E)$ we have
$$ \sup_{S\in \cS^\infty} \E\Big\n \sumk\ga_k S_k Th_k\Big\n^2 <\infty,$$
This induces a well-defined linear operator 
$$ U: \ga(H,E)\to l^\infty(\cS^\infty; L^2(\O;E)),$$
which is bounded by the closed graph theorem. This means that for some
constant $C\geq 0$ we have
$$
\sup_{S\in \cS^\infty} \E\Big\n \sumk \ga_k S_k Th_k\Big\n^2 
\leq C^2\n T\n_{\ga(H,E)}^2.
$$
Now fix  arbitrary $S_1,\dots,S_n\in\cS$ 
and $x_1,\dots,x_n\in E$, and define
$T\in \ga(H,E)$ by $Th_k = x_k$ for $k=1,\dots,n$ and $Th_k=0$ for
$k\geq n+1$. Choosing $S_k\in\cS$ for $k\geq n+1$ arbitrary, we obtain
$$ 
\bal
\E\Big\n \sum_{k=1}^n \ga_k S_k x_k\Big\n^2 
& \leq \E\Big\n \sumk \ga_k S_k Th_k\Big\n^2 
 \leq C^2 \n T\n_{\ga(H,E)}^2 
= C^2 \E\Big\n \sum_{k=1}^n \ga_k x_k\Big\n^2.
\eal
$$
\end{proof}

As an immediate consequence of the previous two results we note:

\begin{corollary}
If $\cT$ is uniformly $\ga$--radonifying in $\B(\tilde H,H)$ 
and $\cS$ is uniformly $\ga$--radonifying in $\B(H,E)$,
then $\cS \cT$ is uniformly $\ga$--radonifying in $\B(\tilde H,E)$ and
$$ \n \cS \cT\n\ugr \leq \n\cS\n\ugr \n\cT\n\ugr.$$
\end{corollary}

If $E$ does not contain a closed subspace isomorphic to $c_0$, then by 
a result of Hoffmann-J\o rgensen and Kwapie\'n
\cite{HJ, Kwapien:c_0}, a Gaussian sum converges in $L^2(\O;E)$ 
if and only if its partial sums are bounded in $L^2(\O;E)$. In
combination with Theorem~\ref{thm:bdd} we obtain the following
equivalent condition for uniform $\ga$--radonification:

\begin{corollary}\label{cor:as}
Let $E$ be a Banach space not containing a copy of $c_0$. 
Then a subset $\cT$ of $\ga(H,E)$ is uniformly $\ga$--radonifying if and
only if there exists a constant $C\geq 0$ such that 
for all integers $n\geq 1$, all orthonormal $h_1,\dots,h_n\in H$, and all
$T_1, \ldots, T_n \in \cT$,
\[
 \E \Bignorm{\sum_{k=1}^n \ga_k T_k h_k}^2 \leq C^2.
\]
In this situation, $\n \cT\n\ugr \leq C$.
\end{corollary}

Here is a simple application. 

\begin{corollary}[Fatou lemma]\label{cor:fatou}
Let $E$ be a Banach space not containing a copy of $c_0$. 
Let $(\cT_n)_{n\geq 1}$ be an increasing sequence of uniformly
$\ga$--radonifying sets in $\ga(H,E)$ satisfying
$ \sup_{n\geq 1} \n \cT_n\n\ugr <\infty.$
Then $\cT := \bigcup_{n\geq 1} \cT_n$
is uniformly $\ga$--radonifying and
$$ \n \cT\n\ugr \leq \sup_{n\geq 1} \n \cT_n\n\ugr.$$
\end{corollary}
\begin{proof}
Let $(T_k)_{k\geq 1}$ be a sequence in $\cT$. For each $m\geq 1$ choose
$N_m\geq 1$ such that $T_1,\dots,T_m\in \cT_{N_m}.$ For 
all orthonormal $h_1,\dots,h_m\in H$ we have
$$ \E\Big\n\sum_{k= 1}^m \ga_k T_k h_k\Big\n^2
\leq \n \cT_{N_m}\n\ugr \leq C,$$
where $C:= \sup_{n\geq 1} \n \cT_n\n\ugr$.
\end{proof}

The condition $c_0\not\subseteq E$ cannot be omitted:

\begin{example}
Let $(e_k)_{k\geq 1}$ and $(u_k)_{k\geq 1}$ denote the standard unit
bases of $\ell^2$ and $c_0$, respectively. It is a classical example of
Linde and Pietsch \cite{LP} that the operator $T\in \B(\ell^2,c_0)$ defined by 
$Te_k = (\ln (k+1))^{-\frac12} u_k$ fails to be $\ga$--radonifying but 
satisfies
$$ \sup_{n\geq 1} \  \E\Big\n\sum_{k= 1}^n \ga_k T e_k\Big\n^2 < \infty.$$
Let $P_k$ be the rank one projection
$e_k\otimes e_k$ in $\ell^2$. Then the sets $\cT_n := \{TP_1,\dots,TP_n\}$ 
satisfy the assumptions of Corollary~\ref{cor:fatou}, but their union
fails to be uniformly $\ga$--radonifying. 
\end{example}

\section{Uniformly $\ga$--radonifying families and compactness in $\ga(H,E)$}
\label{sec:compact}

For an operator $T\in \ga(H,E)$ we define $\mu_T$ as the distribution
of the random variable $\sumk \ga_k Th_k,$ where $(h_k)_{k\geq 1}$ is
an arbitrary orthonormal basis of $H$. The measure $\mu_T$ is a 
centred Gaussian Radon measure on $E$ which does not depend
on the choice of the basis $(h_k)_{k\geq 1}$ and whose covariance
operator equals $TT\s$. For more information on Gaussian measures we
refer to \cite[Chapter 3]{Bogachev:gaussian}, whose terminology we
follow. 

The first result of this section gives a necessary and sufficient
condition for relative compactness in the space $\ga(H,E)$. In a rephrasing in terms of sequential convergence in $\ga(H,E)$, this result is due to Neidhardt \cite{Neidhardt:thesis}. For reasons of self-containedness we shall give a different proof
based on a characterisation of compactness in Lebesgue-Bochner spaces.

\begin{theorem}\label{thm:g-rad-conv}
For a subset $\cT$ of $\ga(H,E)$ the following assertions are equivalent:
\begin{enumerate}
\item The set $\cT$ is relatively compact in $\ga(H,E)$;
\item The set $\{\mu_T\suchthat  T\in\cT\}$ is uniformly tight, and for all $x\s\in
E\s$ the set $\{T\s x\s\suchthat  T\in\cT\}$ is relatively compact in $H$.
\end{enumerate}
\end{theorem}

\begin{proof}
Let $(h_k)_{k\geq 1}$ be a fixed orthonormal basis of $H$
and define, for $T\in \cT$, the random variable $X_T\in L^2(\O;E)$ as
$X_T:= \sumk \ga_k Th_k$. Since $T\mapsto X_T$ defines an isometry
from  $\ga(H,E)$ onto a closed subspace of $L^2(\O;E)$, $\cT$ is relatively
(weakly) compact in $\ga(H,E)$ if and only if $\{X_T\suchthat  T\in\cT\}$ is 
relatively (weakly) compact in $L^2(\O;E)$. 
With this in mind, the proof of the theorem will be 
based on the following compactness result of Diaz and Mayoral
\cite{DM}: for $1\leq p<\infty$,  a subset $A$ of $L^p(\O;E)$ is
relatively weakly compact if and only if the following three
conditions are satisfied:

\ben
\item[(i)] The set $A$ is uniformly $p$-integrable;
\item[(ii)] The set of distributions $\{\mu_f\suchthat\, f\in A\}$ is uniformly tight;
\item[(iii)] The set $\{\lb f,x\s\rb\suchthat f\in A\}$ is relatively weakly compact
in $L^p(\O)$ for all $x\s\in E\s$.
\een
An elementary proof of this result, valid for arbitrary 
Banach function spaces over $(\O,\P)$ with order continuous norm, may be found in \cite{Ne-BFS}.

(a)$\Rightarrow$(b): The uniform tightness of  
the set $\{\mu_T\suchthat T\in\cT\}$ follows from the Diaz-Mayoral result. 
For every $x\s\in E\s$ the set
$\{T\s x\s\suchthat T\in\cT\}$ is relatively compact in $H$, since it is the 
image of the relatively compact set $\cT$ under the continuous mapping
from $\ga(H,E)$ into $H$ given by $T\mapsto T\s x\s$.

(b)$\Rightarrow$(a): 
The relative compactness of $\{T\s x\s\suchthat T\in\cT\}$ in $H$
implies the relative compactness in $L^2(\O)$ 
of the random variables $\{\lb X_T,x\s\rb\suchthat  T\in \cT\}$. To see this,
just note that
$$ \n \lb X_{T_1},x\s\rb - \lb X_{T_2},x\s\rb\n_{L^2(\O)}^2
=\n \lb X_{T_1-T_2},x\s\rb\n_{L^2(\O)}^2
= \n T_1\s x\s - T_2\s x\s\n_H^2.
$$
By \cite[Theorem 3.8.11]{Bogachev:gaussian}, uniformly tight families
of centred Gaussian $E$-valued random variables are uniformly square
integrable and therefore (b) implies (a) by another application of the
compactness result of Diaz and Mayoral. 
\end{proof}

 \begin{corollary}\label{cor:dom}
 Let $\cT$ be a subset in $\ga(H,E)$ which is {\em dominated} by some fixed
 element $S\in \ga(H,E)$, in the sense that for all $x\s\in E\s$,
 $$ 0\leq \n T\s x\s\n_H \leq \n S\s x\s\n_H.$$
 Then the following assertions are equivalent:
\begin{enumerate}
\item $\cT$ is relatively compact in $\ga(H,E)$; 
\item For all $x\s\in E\s$ the set 
$\{T\s x\s\suchthat T\in\cT\}$ is relatively compact in $H$.
\end{enumerate} 
\end{corollary}
\begin{proof} A standard domination result 
 for Gaussian measures (see \cite[Theorem 8.8]{NeeCanberra}) implies that if $\cT$ is dominated,
 then the family $\{\mu_T\suchthat  T\in \cT\}$
 is uniformly tight. The result now follows from Theorem~\ref{thm:g-rad-conv}. 
 \end{proof}

The second main result of this section is the following characterisation of 
relative compactness in $\ga(H,E)$ of uniformly $\ga$--radonifying families.

\begin{theorem}\label{thm:unif-gamma-compact}
Let $\cT$ be uniformly
$\ga$--radonifying subset of $\ga(H,E)$. The following assertions hold:
\begin{enumerate}
\item\label{item:unif-gamma-strongly-compact} 
      $\cT$ is relatively compact in $\ga(H, E)$ if and only if
      $\cT h$ is relatively compact in $E$ for all $h\in H$;
\item\label{item:unif-gamma-compact-wot-compact} 
      $\cT$ is relatively weakly compact in $\ga(H,E)$  if and only if
      $\cT h$ is relatively weakly compact in $E$ for all $h\in H$.
\end{enumerate}
\end{theorem}
\begin{proof}
The relative (weak) compactness of $\cT$ in
$\ga(H,E)$ clearly implies the relative (weak) compactness of $\cT h$ in $E$
for all $h\in H$, so we only need to prove the converse statements.
Throughout the proof we fix an orthonormal basis $(h_k)_{k\geq 1}$ of
$H$. 

The proof of~\ref{item:unif-gamma-strongly-compact} is divided into two steps.

\smallskip
{\em Step 1} -- Fix $n\geq 1$ and let $P_n$ be the orthogonal projection in
$H$ onto the span $H_n$ of $h_1,\ldots,h_n$. 
The set $\cT_n = \{TP_n \suchthat  T\in\cT\}$ is
relatively compact in $\ga(H,E)$ by Corollary~\ref{cor:dom}.

\smallskip
{\em Step 2} -- Assume that $\cT$ is not relatively compact; we shall
prove that $\cT$ is not uniformly $\ga$--radonifying. 

Since $\cT$ is not totally bounded, we can find an $\eps>0$ such that
$\cT$ cannot be covered with finitely many $3\eps$-balls. 
We shall construct an increasing sequence of positive integers 
$0= N_0 < N_1 < \ldots$ and a sequence $T_1,T_2,\ldots$ of elements of
$\cT$ such that for all $m\geq 1$ we have
\[  
  \E \Bignorm{\sum_{m=N_{k-1}+1}^{N_k} \ga_k T_m h_k}^2 \geq 4\eps^2.
\]

The $3\eps$-ball with centre $0$ does not cover $\cT$, and therefore we
may pick $T_1\in\cT$ such that $\norm{T_1}_{\ga(H,E)}\geq 3\eps$. Choose
the index $N_1\geq 1$ in such a way that 
\[ 
  \E \Bignorm{\sum_{k=1}^{N_1} \ga_k T_1 h_k}^2 \geq \eps^2.
\]
We claim that for some $T_2\in \cT$ we have
\[
  \E \Bignorm{\sum_{k=N_1+1}^\infty \ga_k T_2 h_k}^2 \geq 4\eps^2.
\]
Suppose this claim was false. Denoting by $Q_{N_1} = I-P_{N_1}$ the
orthogonal projection onto $H_{N_1}^\perp$, this would mean that
$ \norm{TQ_{N_1}}_{\ga(H,E)} < 2\eps$ for all $T\in \cT$.
Then for all $T\in\cT$ we have 
$$ T = TP_{N_1}+ TQ_{N_1} \in \cT_{N_1} + \B(2\eps),$$ 
where $\B(2\eps)$ is the $2\eps$-ball in $\ga(H,E)$ centred at $0$. By Step
1 we can cover $\cT_{N_1}$ with finitely many $\eps$-balls, and
therefore   we can cover $\cT$ with finitely many $3\eps$-balls. This
contradiction proves the claim. Now choose the index $N_2 \geq N_1+1$
in such a way that 
\[ 
   \E \Bignorm{\sum_{k=N_1+1}^{N_2} \ga_k T_2 h_k}^2 \geq \eps^2.
\]
It is clear that this construction can be continued inductively.

Let $S_k := T_m$ if $N_{m-1}+1\leq k\leq N_m$ for some $m\geq 1$.
Then $(S_k)_{k\geq 1}$ is a sequence in $\cT$ for which the sum
$\sum_k \ga_k S_k h_k$ fails to converge. 
\smallskip

Next we prove~\ref{item:unif-gamma-compact-wot-compact}. 
We say that the sequence
$(y_n)_{n\geq 1}$ is a {\em convex tail subsequence} of a sequence 
$(x_n)_{n\geq  1}$ in $E$ if each $y_n$ is a convex combination of elements of the tail
sequence $(x_k)_{k\geq n}$. Note that if $\limn x_n = x$ strongly or
weakly, then also $\limn y_n = x$ strongly or weakly.
We shall use of the following weak compactness criterium 
\cite[Corollary 2.2]{DiestelRuessSchachermayer:compactness}: a subset
$K$ of a Banach space $X$ is relatively weakly compact if and only if
every sequence in $K$ has a strongly convergent convex tail
subsequence.

After these preparations we turn to the proof of
\ref{item:unif-gamma-compact-wot-compact}. Let $(T_k)_{k\geq 1}$ be a
sequence in $\cT$.  By a diagonal argument we find a subsequence
$(T_{k_j})_{j\geq 1}$ such that the weak limit $\limj T_{k_j} h_n$
exists for every $n\geq 1$. By a standard corollary to the Hahn-Banach
theorem and a diagonal argument we find a convex tail subsequence
$(S_j)_{j\geq 1}$  of $(T_{k_j})_{j\geq 1}$ such that the strong limit 
$\limj S_j h_n $ exists for every  $n\geq 1$. By the uniform
boundedness of $\cT$, the strong limit $Sh := \limj S_{j} h$  exists
for all $h\in H$. Now part (a) implies that $\limj S_{k_j} = S$ in
$\ga(H,E)$. Hence, by the above criterium, $\cT$ is relatively weakly
compact. 
\end{proof}

The following example shows that uniformly $\ga$--radonifying families in
$\ga(H,E)$ need not
be relatively compact in $\ga(H,E)$, even in the case where $E$ is a Hilbert space.

\begin{example}  
Let $(e_k)_{k\geq 1}$ be the standard unit basis of $\ell^2$ and fix an
arbitrary nonzero element $h\in \ell^2$. We check that the family
$$ \cT := \{h\otimes e_j\suchthat  j\geq 1\}$$ is a uniformly
$\ga$--radonifying subset of $\ga(\ell^2):= \ga(\ell^2,\ell^2)$. Taking this for granted for
the moment, noting that $\{Th\suchthat  T\in \cT\}$ fails to be
relatively compact in $\ell^2$  it follows from Theorem
\ref{thm:unif-gamma-compact} that $\cT$ fails to be relatively compact
in $\ga(\ell^2)$.  

If $(T_k)_{k\geq 1}$ is a sequence in $\cT$, say $T_k = h\otimes e_{j_k}$, 
then for all $1\leq M\leq N$ we have 
$$
\bal
 \E \Big\n \sum_{k=M}^N \ga_k T_k e_k\Big\n^2
 = \E \Big\n \sum_{k=M}^N \ga_k (h,e_k) e_{j_k}\Big\n^2
 = \sum_{k=M}^N  \n (h,e_k) e_{j_k}\n^2
= \sum_{k=M}^N |(h,e_k)|^2.
\eal 
$$
As $M,N\to \infty$ the right-hand side tends to $0$, which proves that
$\sumk\ga_k T_k e_k$ converges in $L^2(\O;\ell^2)$.
\end{example} 

Our next aim is to show that for every Banach space $E$
there exists a relatively compact $\cT$ set in $\ga(\ell^2,E)$ 
which fails to be uniformly $\ga$--radonifying. 

\begin{example}\label{ex:not-ugr}
Let $(h_k)_{k\geq  1}$ denote the standard unit basis of $\ell^2$. 
Define $S\in \B(\ell^2)$ to be the right shift, i.e. $S h_k = h_{k+1}$
for all $k\geq 1$. For $T\in\ga(\ell^2,E)$ let $\cS_T:=
\{TS^n\suchthat  n\in\N\}$. This set  is bounded in $\ga(\ell^2,E)$
and for all $n\in\N$ we have
\[
  \norm{S^{n*} T\s x\s}_{\ell^2}^2 \leq
  \norm{T\s x\s}_{\ell^2}^2.
\]
Also, for all $x\s\in E\s$ we have $\limn S^{n*} T\s x\s = 0$ strongly,
and therefore  $\cS_T$ is relatively compact by Corollary~\ref{cor:dom}.

In what follows we take for $E$ the scalar field ${\mathbb K}$.
Define $M_1 = 1$ and, inductively, $M_{n+1} :=  M_n + n$ for $n\geq 1$.
Consider the operators $T_n:\ell^2\to {\mathbb K}$ defined by
$T_n h_{M_{n+1}} = 1$ and $T_n h_k =0$ for $k\neq M_{n+1}$. Trivially, 
this operator is
$\ga$--radonifying with  $ \norm{T_n }_{\ga(\ell^2,{\mathbb K})} = 1$.
Let $m_k^n:= n - j$ for $k=M_n+j$, $j=1,\dots,n$. Then,
\[
  \Bigl(\E \Bignorm{\sum_{k\geq 1} \ga_k T_n S^{m_k^n} h_k}^2\Bigr)^\einhalb
= \Bigl(\E \Bignorm{\sum_{k=M_n+1}^{M_{n+1}} \ga_k T_n
     h_{M_{n+1}}}^2\Bigr)^\einhalb 
= (M_{n+1}-M_{n})^\einhalb = n^\einhalb
\]
and similarly,
$$
  \Bigl(\E \Bignorm{\sum_{k\geq 1} \ga_k T_{n'} S^{m_k^n}
    h_k}^2\Bigr)^\einhalb = 0, \qquad n\not=n'.
$$
Define $T:\ell^2\to {\mathbb K}$ by 
$T:=  \sum_{n\geq 1} \frac1{n^2} T_{2^n}.$
Note that $T\in\ga(\ell^2,{\mathbb K})$.
By the contraction principle, for all $n\geq 1$ we have 
$$\bal
      \Bigl(\E \Bignorm{\sum_{k=1}^{M_{2^n+1}} \ga_{k} T S^{m_k^{2^n}}
         h_k}^2\Bigr)^\einhalb
&\geq \Bigl(\E \Bignorm{\sum_{k=M_{2^n}+1}^{M_{2^n+1}} \ga_k T
         S^{m_k^{2^n}} h_k}^2\Bigr)^\einhalb \\ 
& =  \frac{1}{n^2} \Bigl(\E \Bignorm{
        \sum_{k=M_{2^n}+1}^{M_{2^n+1}} \ga_{k} T_{2^n}S^{m_k^{2^n}}
        h_k}^2\Bigr)^\einhalb
 =    \frac1{n^{2}}\cdot 2^{\frac{n}2}.
\eal
$$ 
Since $n$ is arbitrary, this implies that the family
$\cS_T = \{TS^m\suchthat  m\in\N\}$
fails to be uniformly $\ga$--radonifying. Note that this family is bounded, hence $\ga$-bounded, in $\calL(\ell^2,\mathbb{K})$.
\end{example}

\section{Laplace transforms}\label{sec:laplace}

Let $I$ be a countable index set. A sequence $(h_i)_{i\in I}$ in a
Hilbert space $H$ is said to be a {\em Hilbert sequence} if there
exists a constant $C>0$ such that for all scalar sequences $\al \in
\ell^2(I)$,
\[
\Bignorm{ \sum_{i\in I} \al_i h_i }_H^2\leq C^2\sum_{i\in I} |\al_i|^2.
\]
The infimum of all admissible constants $C$ will be called the 
{\em Hilbert constant} of the  sequence $(h_i)_{i\in I}$, 
cf. \cite[Section 1.8]{Young:nonharmonic}.
The usefulness of this notion is explained by the following result
\cite[Proposition 2.1]{HNV}:

\begin{proposition}\label{prop:HNV} Let $T\in \ga(H,E)$ be given.
If $(h_i)_{i\in I}$ is a Hilbert sequence in $H$, then the Gaussian sum
$\sum_{i\in I} \ga_i Th_i $ converges in $L^2(\O;E)$ and we have
$$ \E\Big\n \sum_{i\in I} \ga_i Th_i\Big\n^2 \leq C^2 \n T\n_{\ga(H,E)}^2,$$
where $C$ is the Hilbert constant of $(h_i)_{i\in I}.$
\end{proposition}

\begin{example}\label{ex:JZ}
Let $(\la_n)_{n\geq 1}$ be a sequence in $\C_+$ which is {\em properly
  spaced} in the sense that 
$$\inf_{m\not=n} \Big|\frac{\la_m - \la_n}{\Re(\la_n)} \Big| > 0.$$ Then the functions 
$$f_n(t):= \sqrt{\Re(\la_n)}
e^{-\la_n t}, \quad n\ge 1,$$ define a Riesz sequence on the closure 
of their span in 
$L^2(\R_+)$, i.e., there are constants $0<c\le C<\infty$ 
such that 
\[c^2\sum_{n\ge 1} |\al_n|^2\le 
\Bignorm{ \sum_{n\ge 1} \al_n f_n }^2\leq C^2\sum_{n\ge 1} |\al_n|^2
\]
for all sequences $(\al_n)_{n\ge 1}\in \ell^2$; see 
\cite[Theorem 1, (3)$\Leftrightarrow$(5)]{JZ}. In particular, $(f_n)_{n\ge 1}$ is a Hilbert sequence in
$L^2(\R_+)$.
From this one easily deduces that for
 any $b>0$ and $\rho\in [0,1)$ the
functions $$f_n(t) = e^{-bt + 2\pi i(n+\rho)t}, \quad n\in\Z,$$
define a Hilbert sequence in $L^2(\R_+)$. This has been shown by direct computation in \cite[Example 2.5]{HNV}, where the bound $1/\sqrt{1-e^{-2b}}$
was obtained for its Hilbert constant. Note that this bound is independent of $\rho$.
\end{example}

The next proposition is well-known and 
shows that a sequence is a Hilbert sequence if it is not `too far' from being orthogonal.
For the reader's convenience we include an elementary proof.

\begin{proposition}\label{prop:working-horse}
Let $(h_n)_{n\in \ZZ}$ be a sequence in $H$. If there exists a function
$\phi: \NN\to \RR_+$ such that for all $n\geq m \in \ZZ$ we have $\bigl| 
(h_n,h_m)_H \bigr| \leq \phi(n-m)$ and $\sum_{j\in\N} \phi(j) <\infty$, then
$(h_n)_{n\in \ZZ}$ is Hilbert sequence. 
\end{proposition}
\begin{proof}
 Let $(\al_n)_{n\in \ZZ}$ be scalars. 
 Then
 \begin{eqnarray*}
   \Bignorm{ \sum_{n=-N}^N \al_n h_n }^2
 &=&   \sum_{n=-N}^N |\al_n|^2 \norm{h_n}^2
       + 2 \Re\sum_{-N\le n<m\le N} \al_n \overline{\al_m}
       (h_n, h_m)_H\\
 &\leq& \phi(0) \sum_{n\in \ZZ} |\al_n|^2
       + 2 \sum_{n<m} |\al_n| |\al_m| \phi(n-m) \\
 &=&   \phi(0) \sum_{n\in \ZZ} |\al_n|^2
       + 2 \sum_{j\geq 1} \phi(j) \sum_{n\in \ZZ} |\al_n|\, |\al_{n+j}| \\
 &\leq& \Big(\phi(0) + 2  \sum_{j\geq 1} \phi(j)\Big) \sum_{n\in \ZZ}
       |\al_n|^2,
 \end{eqnarray*}
 where the last estimate follows from the Cauchy-Schwarz inequality.
 \end{proof}

As a special case we have the following example, which will be needed
in the proof of Theorem~\ref{thm:sqrt-la2}.

% %%%%%%%%%%%%%%%%%
\begin{example}\label{ex:hilbert-sequence-laplace}
% %%%%%%%%%%%%%%%%% 
Let $\al \in (0,\einhalb]$, $r>0$, and $\vartheta \in (-\frac{\pi}{2},
\frac{\pi}{2})$. Let $\mu_n = r 2^{n} e^{i \vartheta}$, $n\in\Z$, and let 
\beq\label{seq}
f_n(s) := \mu_n^{\al} e^{-\mu_n s}, \quad s\in\R_+.
\eeq Then $(f_n)_{n\in \ZZ}$ is
a  Hilbert sequence in $L^2(\RR_+)$.
Indeed, for $n,m \in \ZZ$, $n\not=m$, we have
\[\bal
 |(f_n, f_m)|
  & \leq \int_0^\infty r^{2\al} 2^{\al(n+m)}
       e^{-r(2^n+2^m)s\cos\vartheta } \, ds
    =  \frac{r^{2\al-1} 2^{\al(n+m)}}{(2^n+2^m)\cos\vartheta}.
\eal\]
Since $\al\in (0,\einhalb]$ we have
\[
    \frac{r^{2\al-1} 2^{\al(n+m)}}{(2^n+2^m)\cos\vartheta}
 \leq \frac{r^{2\al-1}}{\cos\vartheta} \frac{2^{\al(n+m)}}{2^{\max(n,m)}}
 \leq \frac{r^{2\al-1}}{\cos\vartheta} 2^{-\al|n-m|},
\]
and Proposition~\ref{prop:working-horse} applies. Notice that for $\al=\einhalb$,
the obtained Hilbert constant estimate is bounded by $C/\cos\vartheta$, where $C$ is a universal constant.
\end{example}

From now on, $H$ is again a separable infinite-dimensional Hilbert space.
The main abstract result of this section reads as follows.

For an operator $\Phi\in \ga(L^2(\R_+;H),E)$ and a function
$f\in L^2(\R_+)$ we define the operator $f(\Phi)\in\ga(H,E)$
by
$$f(\Phi)h := \Phi(f\otimes h), \quad h\in H.$$
Below we shall apply this definition to the functions $f(t)= e^{-\la t}$ with $\Re\la>0$ to in order to define the `Laplace transform' of $\Phi$.

\begin{theorem}\label{thm:laplace-abstrakt}
Let $(f_i)_{i\in I}$ be a Hilbert sequence in
$L^2(\RR_+)$ with Hilbert constant $C$. 
Then for all $\Phi \in \ga(L^2(\RR_+;H), E)$ the
family $$\cT := \{f_i(\Phi) \suchthat i\in I\}$$ is
uniformy $\ga$--radonifying
and we have
$$ \n \cT\n\ugr \leq C \n \Phi\n_{\ga(L^2(\R_+;H),E)}.$$
\end{theorem} 
\begin{proof} Fix an
orthonormal basis $(h_k)_{k\geq 1}$ in $H$ and 
let $(i_k)_{k\geq 1}$ be an arbitrary sequence in $I$. 
Put $J := \{i\in I: \ i_k =i$ for some $k\in K\}$. For each $i\in J$,
put $K(i) := \{k \ge 1 \suchthat i_k = i \}$.
Fix a Gaussian sequence $(\ga_k)_{k\ge 1}$
on a probability space $(\Om,\P)$, 
as well as a doubly indexed Gaussian sequence $(\ga_{ik}')_{i\in I,\,k\ge 1}$
on another probability space $(\Om',\P')$.
We have
$$ \bal
   \EE \Bigl\n\sumk \ga_k f_{i_k}(\Phi)h_k\Bigr\n^2
& =   \EE'\Bigl\n\sumk \ga_{i_kk}' f_{i_k}(\Phi) h_k\Bigr\n^2
\\ &=  \EE' \Bigl\n \sum_{i\in J} \sum_{k \in K(i)} \ga_{ik}'
       {f_i}(\Phi) h_k \Bigr\n ^2 
 \leq  \EE' \Bigl\n  
     \sum_{i \in I} \sumk\ga_{ik}'{f_i}(\Phi) h_k\Bigr\n^2.
\eal $$
To prove convergence of the double sum on the right-hand side we note that 
that the sequence
$(f_i\otimes h_k)_{i\in I,\,k\geq 1}$ is a
Hilbert sequence in $L^2(\RR_+;H)$ with Hilbert constant $C$. Indeed, this follows from
$$
\bal
  \Bignorm{\sum_{i\in I} \sumk \al_{ik}
    f_i \otimes h_k}_{L^2(\RR_+;H)}^2
&= \int_0^\infty \Bignorm{ \sum_{i\in I}\sumk \al_{ik} f_i(t) h_k}_H^2 \,dt\\
&= \int_0^\infty \sumk\Big|\sum_{i\in I}  \al_{ik} f_i(t)\Big|^2 \,dt
\\ & = \sumk\Big\n\sum_{i\in I}  \al_{ik} f_i\Big\n_{L^2(\R_+)}^2
\\ & \leq C^2\sum_{i\in I} \sumk |\al_{ik}|^2.
\eal
$$
Hence by Proposition~\ref{prop:HNV},
$$ \E\Big\n \sumk \ga_k f_{i_k}(\Phi) h_k\Big\n^2
\leq  
C^2 \n \Phi\n_{\ga(L^2(\R_+;H),E)}^2 .$$
\end{proof}
 
We shall present three applications of this result.

The {\em Laplace transform} of an operator $\Phi \in \ga(L^2(\RR_+;H), E)$ 
is the function $\widehat\Phi:\C_+\to \ga(H,E)$ defined by
$$ \widehat\Phi(\la)h:= e_\la(\Phi)h = 
\Phi(e_\la\otimes h),\quad h\in H,$$
where $\C_+ :=  \{\Re\la>0\}$ and $e_\la(t):= e^{-\la t}$ for $t\in\R_+$ and
$\la\in\C_+$.
 
An operator $\Phi\in \ga(L^2(\RR_+;H), E)$ is said to be {\em $H$-strongly $L^1$-representable}
if for all $h\in H$ there exists a function $\phi_h\in L^1(\R_+;E)\cap L^2(\R_+;E)$
such that for all $f\in L^2(\R_+)$ we have
$$\Phi (f\otimes h) = \int_0^\infty f(t) \phi_h(t)\,dt.$$ 
Under this assumption we have
$$\widehat{\Phi}(\la)h = \widehat{ \phi_h}(\la), \quad \Re\la> 0.$$ 

\begin{theorem}[$\ga$--Riemann-Lebesgue lemma]\label{thm:g-RL}
If $\Phi\in \ga(L^2(\RR_+;H), E)$ is $H$-strong\-ly $L^1$-representable,
then $$\limn\n\widehat \Phi(\la_n)\n_{\ga(H,E)} = 0$$
for any sequence $(\la_n)_{n\geq 1}$ in $\C_+$ such that 
$(e_{\la_n})_{n\geq 1}$ is a Hilbert sequence.
\end{theorem}
\begin{proof}
Let $(\la_n)_{n\geq 1}$ be a sequence in $\C_+$ as stated. 
By \Fref{thm:laplace-abstrakt}, 
the family  $\{ \widehat \Phi(\la_n)\suchthat n\geq 1\}$
is uniformly $\ga$--radonifying. 
Moreover, for all $h\in H$
we have  $$\lim_{|\la|\to\infty} \widehat \Phi(\la) h
=\lim_{|\la|\to\infty} \widehat{ \phi_h}(\la)=0$$
by the Riemann-Lebesgue lemma. Consequently, for every $h\in H$ the set 
$\{\widehat \Phi(\la)h\suchthat \la\in \C_+\}$ is relatively compact in
$E$. \Fref{thm:unif-gamma-compact} then shows that 
$\{\widehat \Phi(\la_n)\suchthat n\geq 1\}$ is relatively compact in 
$\ga(H, E)$. Therefore, $ \lim_{n\to\infty}\widehat \Phi(\la_n)h= 0$
implies $ \lim_{n\to\infty}\widehat \Phi(\la_n)= 0$ in $\ga(H, E)$.
\end{proof}

In particular we obtain that 
if $\Phi\in \ga(L^2(\RR_+;H), E)$ is $H$-strong\-ly $L^1$-represent\-able,
then  for all $b>0$ we have
$$\lim_{|s|\to \infty}\n\widehat \Phi(b+is)\n_{\ga(H,E)} = 0.$$

In the next two applications of Theorem~\ref{thm:laplace-abstrakt} we consider the uniform $\ga$--radonif\-ication 
of Laplace transforms in right half-planes and sectors, respectively.

Let $S = \{\lambda\in\C\suchthat  0<\Re\lambda<1\}$. If $N:\overline
S\to \calL(E,F)$ is strongly continuous and bounded on $\overline S$
and harmonic on $S$,  then by the Poisson formula for the strip 
\cite{Widder}, cf. also \cite{vanNeervenWeis:asymptotic-scp},
we have, for $\lambda=\al+i\beta$ with $0<\al<1$ and
$\beta\in\R$,
\beq \label{eq:poisson} 
N(\lambda)x = \sum_{j=0,1} \int_{-\infty}^\infty 
 P_j(\al, \beta-t)N(j+it)x\,dt, \qquad x\in E,
\eeq
where
$$ P_j(\al,s) =\frac1\pi \frac{e^{\pi s} \sin(\pi\al)}{\sin^2(\pi\al) +
(\cos(\pi\al)-(-1)^je^{\pi s})^2}.$$

\begin{theorem}[Uniform $\ga$--radonification in half-planes]
\label{thm:sqrt-la}
Let $\Phi \in \ga(L^2(\RR_+;H), E)$ be given. For all $b>0$ 
the family $$\cT_b^\Phi:= \bigl\{\widehat \Phi (\la) \suchthat \Re\la \geq
b\bigr\}$$
is uniformly $\ga$--radonifying in $\ga(H, E)$ and
$$ \n \cT_b^\Phi\n\ugr \leq \frac{C}{\sqrt{b}}\norm{\Phi}_{\ga(L^2(\RR_+;H), E)}, $$ 
where $C$ is a universal constant.
\end{theorem}

\begin{proof}
By  Example~\ref{ex:JZ} and a substitution (cf. \cite[Theorem 3.1]{HNV}),
for $\si \in [\frac{1}{2}b, \frac32 b]$ and $\rho\in [0,1)$ fixed,
the sequence $(g_n)_{n\in \ZZ}$ given by
\[
g_n(t) := e^{-\si t + i(n{+}\rho)b t}, \quad t\in\R_+,
\]
is a Hilbert sequence with Hilbert constant $\leq
\tfrac{C}{\sqrt{b}}$,  where $C = \sqrt{\tfrac{2\pi
  e^{2\pi}}{e^{2\pi}-1}}$.
Consequently, \Fref{thm:laplace-abstrakt} shows the uniform
$\ga$--radonification of the set $\{ \widehat
\Phi(\si+i(n{+}\rho)b)\suchthat n\in \ZZ\}$ with constant 
$\leq\tfrac{C}{\sqrt{b}}\norm{\Phi}_{\ga(L^2(\RR_+;H), E)}$.

Let $(h_k)_{k\geq 1}$ be an orthonormal 
basis of $H$ and let $(\la_k)_{k\geq 1}$ be a sequence on the line 
$\{\Re\la = b\}$, say $\la_k = b +i(n_k+\rho_k)b$ with $n_k\in \Z$ and
$0\leq \rho_k<1$.

Fix indices $1\leq M\leq N$. Following the argument of \cite[Theorem
4.3]{vanNeervenWeis:asymptotic-scp}, the 
Poisson integral formula \eqref{eq:poisson} can be used with
$N(\la) = \widehat\Phi((\tfrac12+\la)b)$
to estimate
$$ 
\begin{aligned}
 \ & \Bigl(\E \Bigl\n\sum_{k=M}^N \ga_k\,\widehat\Phi(\lambda_k) h_k\Bigr\n\Bigr)^\frac12
\\ & \ \ =  \Bigl\n\sum_{j=0,1}\sum_{k=M}^N \ga_k\, \int_{-\infty}^\infty 
P_j(\tfrac12, n_k+\rho_k-t) \widehat\Phi((\tfrac12 +j)b+itb)h_k\,dt\Bigr\n_{L^2(\O;E)}
\\ & \ \ \leq \sum_{j=0,1}  \int_{-\infty}^\infty \Bigl\n
\sum_{k=M}^N \ga_k\,
  P_j(\tfrac12, \rho_k -\tau ) \widehat\Phi((\tfrac12 +j)b+i(n_k+\tau)b)h_k\Bigr\n_{L^2(\O;E)}\,d\tau
\\ & \ \ \leq \sum_{j=0,1}  \int_{-\infty}^\infty \sup_{\rho\in [0,1)} 
 P_j(\tfrac12, \rho -\tau) \, \Bigl\n\sum_{k=M}^N \ga_k\, 
 \widehat\Phi((\tfrac12+j)b+i(n_k+\tau)b)h_k\Bigr\n_{L^2(\O;E)}\,d\tau.
\end{aligned}
$$
In the last estimate we used the contraction principle. 
For fixed $\tau\in\R$ we have
$$\lim_{M,N\to\infty}  \Bigl\n\sum_{k=M}^N \ga_k\, 
 \widehat\Phi((\tfrac12+j)b+i(n_k+\tau)b)h_k\Bigr\n_{L^2(\O;E)} =0$$
since $\big\{ \widehat\Phi((\tfrac12+j)b+i(n+\tau)b): \ n\in\Z\big\}$
is uniformly $\ga$--radonifying, and
$$
\bal 
\ &  \sup_{\rho\in [0,1)} 
 P_j(\tfrac12, \rho -\tau) \, \Bigl\n\sum_{k=M}^N \ga_k\, 
 \widehat\Phi((\tfrac12+j)b+i(n_k+\tau)b)h_k\Bigr\n_{L^2(\O;E)}
\\ & \hskip2cm \leq 
\frac{C}{\sqrt{b}} \Big(\sup_{\rho\in [0,1)} 
 P_j(\tfrac12, \rho -\tau) \Big)\norm{\Phi}_{\ga(L^2(\RR_+;H), E)}.
\eal $$
Since the right-hand side is an integrable function of $\tau$ we may
apply dominated convergence to conclude that
$$
\lim_{M,N\to\infty}\Bigl(\E \Bigl\n\sum_{k=M}^N
\ga_k\,\widehat\Phi(\lambda_k) h_k\Bigr\n\Bigr)^\frac12 =0.$$
This shows that
$\{\widehat \Phi (\la)\suchthat  \Re\la = b\}$ is uniformly
$\ga$--radonifying. Moreover, taking $M=1$ and letting $N\to\infty$ in
the above estimates, we obtain the bound
$$
\Bigl(\E \Bigl\n\sumk \ga_k\,\widehat\Phi(\lambda_k) h_k\Bigr\n\Bigr)^\frac12
\leq 2\sup_{j=0,1}
\frac{C}{\sqrt{b}}
\Big(\int_{-\infty}^\infty \sup_{\rho\in [0,1)} P_j(\tfrac12, \rho-\tau)\,d\tau \Big)
\norm{\Phi}_{\ga(L^2(\RR_+;H), E)}.
$$
This proves that $\{\widehat{\Phi}(\la)\suchthat  \Re \la = b\}$ is
uniformly $\ga$--radonifying with constant $\leq
C'/{\sqrt{b}}\norm{\Phi}_{\ga(L^2(\RR_+;H), E)}$, where $C'$ is
universal. By \Fref{prop:jordan}, $\{\widehat{\Phi}(\la)\suchthat 
\Re \la \geq b\}$ is then uniformly $\ga$--radonifying with at most
twice this constant.  
\end{proof}

Combining this theorem with Corollary~\ref{cor:gbdd2} we recover 
\cite[Theorem 3.4]{vanNeervenWeis:asymptotic-scp}, which asserts
that the Laplace transform of $\Phi$ is $R$--bounded on $\{\Re\la\geq
b\}$ for all $b>0$, with an $R$--bound  of order
$O(\frac{1}{\sqrt{b}})$ as $b\downarrow 0$.  In view of Example
\ref{ex:not-ugr}, Theorem~\ref{thm:sqrt-la} represents a genuine
strengthening of this result.

\medskip
Next we turn to the uniform $\ga$--radonification of Laplace transforms in sectors. 
Before we can state and prove our main result in this direction, Theorem~\ref{thm:sqrt-la2},
we introduce some notations.

For $0<\vartheta<\pi$ and $0<r<R $ we define 
$$S_{\vartheta} : = \{z\in \C\suchthat  z\not=0, \
|\arg z|< \vartheta\},$$ where the argument is taken in $(-\pi,\pi)$.

% %%%%%%%%%%%%%%%%%
\begin{theorem}[Uniform $\ga$--radonification in sectors]\label{thm:sqrt-la2}
% %%%%%%%%%%%%%%%%%
Let $\Phi\in \ga(L^2(\RR_+;H), E)$ be given. 
For all $0<\vartheta<\pihalbe$
the family 
$$\cT_\vartheta^\Phi := \{ \sqrt{\la}\, \widehat \Phi(\la)\suchthat \la \in
S_\vartheta\}$$ is uniformly $\ga$--radonifying in $\ga(H,E)$ and
$$ \n \cT_\vartheta^\Phi\n\ugr \leq \frac{C'}{\cos\vartheta}\n
\Phi\n_{\ga(L^2(\R_+;H),E)},$$
where $C'$ is a universal constant.
\end{theorem}
\begin{proof}
The proof follows the lines of Theorem~\ref{thm:sqrt-la}, the difference being that instead of using Example~\ref{ex:JZ} we now use Example~\ref{ex:hilbert-sequence-laplace}.

Fix $\vartheta<\theta<\pihalbe$ such that $\cos\theta\ge \frac12\cos\vartheta$.
One obtains that for any $r>0$ fixed,
the sequences $(f_n^+)_{n\in \ZZ}$ and $(f_n^-)_{n\in \ZZ}$ given by
$$f_n^\pm(t) = \sqrt{\mu_{n}^\pm} e^{-\mu_{n}^\pm t}$$
with $$ \mu_{n}^\pm = r 2^{n} e^{\pm i \theta}$$
are Hilbert sequences whose Hilbert constants are bounded by $C/\cos \theta$,
where $C$ is a universal constant.
Hence, arguing along the lines of Theorem~\ref{thm:laplace-abstrakt},
we obtain that the sequence $(\sqrt{\mu_n}\widehat\Phi(\mu_n))_{n\in \Z}$
is uniformly $\ga$--radonifying, with bound $C/\cos \theta$.

By a Poisson transform argument (e.g., by using the logarithm to conformally map sectors to strips and then using the argument of Theorem~\ref{thm:sqrt-la}),
we obtain that $\sqrt{\la}\widehat\Phi(\la)$
is uniformly $\ga$--radonifying on the sector $S_{\vartheta}$, with a bound $C'/\cos\theta \le   2C'/\cos \vartheta$, where $C'$ is another universal constant. 
\end{proof}

\section{The stochastic Weiss conjecture}\label{sec:weiss}

Let $A$ be the generator of a $C_0$--semigroup $S = (S(t))_{t\geq 0}$
on a Banach space $E$ and let $s(A)$, $s_0(A)$, and $\om_0(A)$ denote 
the spectral bound, the abscissa of uniform boundedness of the resolvent, 
and growth bound of $A$, respectively:
$$
\bal 
    s(A) &= \sup\{\Re\la\suchthat  \la\in\sigma(A)\},\\
  s_0(A) &= \inf\{\om>\sigma(A)\suchthat  \sup_{{\rm Re}\la\geq \om}\n
              R(\la,A)\n < \infty\},\\
\om_0(A) &= \inf\Big\{\om\in\R\suchthat  \n S(t)\n\leq
             Me^{\om t}\hbox{ for some $M\geq 1$ and all $t\geq 0$}\Big\}.
\eal
$$ 
Here, $R(\la,A) := (\la-A)^{-1}$. Recall that 
$-\infty\leq s(A)\leq s_0(A)\leq \om_0(A)<\infty$.

It is shown in \cite{vanNeervenWeis:stoch-int} that 
the linear stochastic Cauchy problem 
$$\left\{
\begin{aligned}
 dU(t) & = AU(t)\,dt + B\,dW_H(t), \qquad t\in[0,T], \\
  U(0) & = u_0,
\end{aligned}
\right.  \leqno{\SCP{A}{B}}
$$
where $(W_H(t))_{t\in[0,T]}$ is an $H$--cylindrical Wiener
process and $B\in\B(H,E)$ is a bounded operator, has a solution if and only if 
for some (all) $t>0$ the $\B(H,E)$-valued function $S(\cdot)B$
{\em represents} an element of $\ga(L^2(0,t;H),E)$,
in the sense that the integral operator
$$ f \mapsto\int_0^t S(s)B f(s)\,ds, \quad f\in L^2(0,t;H),$$
belongs to $\ga(L^2(0,t;H),E)$.
In this situation the solution is unique up to modification. For
the precise notion of `solution'  as well as other
unexplained terminology we  refer to \cite{vanNeervenWeis:stoch-int}.

As an application of Theorem~\ref{thm:g-RL}
we obtain the following necessary condition for the existence of 
solutions to the problem $\SCP{A}{B}$.

\begin{theorem} \label{thm:sol}
If the problem $\SCP{A}{B}$
has a solution, then for all $b>s_0(A)$  we have
$$ \lim_{|s|\to\infty} \n R(b+is,A)B\n_{\ga(H,E)} = 0.$$
\end{theorem}
\begin{proof}
By \cite[Proposition 4.5]{vanNeervenWeis:asymptotic-scp}, for $b>\om_0(A)$ the
integrable $\B(H,E)$-valued function $t\mapsto 
e^{-bt} S(t)B$ represents an element of $\ga(L^2(\R_+;H),E)$.
Hence
for $b>\om_0(A)$ the assertion is an immediate consequence of
Theorem~\ref{thm:g-RL}. For $b> s_0(A)$ 
the result then follows by a standard resolvent identity argument.
\end{proof}

Note that 
we did not assume that $B\in\ga(H,E)$. Indeed, 
in many examples the problem $\SCP{A}{B}$ admits a solution
without such an assumption on $B$. For operators $B\in\ga(H,E)$ 
the theorem is trivial, since then we 
may apply the Riemann-Lebesgue lemma in the space $L^1(\R_+;\ga(H,E))$.

We recall the fact, proved in \cite[Proposition
4.4]{vanNeervenWeis:asymptotic-scp}, that the problem $\SCP{A}{B}$ admits an invariant measure if and only if
the function $t\mapsto S(t)B$ represents an element of $\ga(L^2(\R_+;H),E)$.
In this situation the mapping $\la\mapsto R(\la,A)B$ extends to an
analytic $\ga(H,E)$-valued function on $\C_+$; this extension is given
by $\la\mapsto \widehat\Phi(\la)$, where $\Phi(t):= S(t)B$. With a
slight abuse of notation we shall write $R(\la,A)B$ for this
extension, keeping in mind that this notation is formal; indeed,
examples can be given where $A$ has spectrum in the open right-half
plane.

As an application of Theorem~\ref{thm:sqrt-la} we obtain the following
necessary conditions for the existence of an invariant measure for the
problem $\SCP{A}{B}$.

\begin{theorem} \label{thm:inv} 
  If the problem $\SCP{A}{B}$ admits an invariant measure, then for
  all $0<\vartheta<\pihalbe$ the family $$\cT_\vartheta :=
  \big\{\sqrt{\la}R(\la,A)B\suchthat \la\in S_\vartheta\big\}$$
  is uniformly $\ga$--radonifying and we have
$$ \n \cT_\vartheta\n\ugr \leq \frac{C_{A,B}}{\cos{\vartheta}},$$
  where $C_{A,B}$ is a constant depending only on $A$ and $B$.
\end{theorem}

We conjecture that the following converse of this theorem holds.

\begin{conjecture}[Stochastic Weiss conjecture]
Let $E$ be a Banach space with finite cotype and assume that the operator $-A$ is injective and sectorial of angle $<\pihalbe$ on $E$ and admits a bounded $H^\infty$--calculus.  The following assertions are equivalent:

\begin{enumerate}
\item the stochastic Cauchy problem
  $\SCP{A}{B}$ admits an invariant measure;
\item The operator $(-A)^{-\frac12}B$ is $\ga$--radonifying;
\item the set
  $\big\{\sqrt{\la}R(\la,A)B\suchthat  \la\in S_\vartheta\}$ is uniformly $\ga$--radonifying for some/all $0< \vartheta<\pihalbe$.
\end{enumerate}
\end{conjecture}

The implication 
(a)$\Rightarrow$(c) follows from Theorem~\ref{thm:inv}, and the implication
(b)$\Rightarrow$(a) can be proved as follows. 
Since $(-A)^{-\frac12}B\in \ga(H,E)$ we may apply \cite[Theorem 6.2]{DNW} 
to obtain that for all $t>0$ the function
$S(\cdot)B = (-A)^\frac12 S(\cdot) ((-A)^{-\frac12}B)$ belongs to
$\ga(L^2(0,t;H),E)$,
with a uniform bound $\sup_{t>0} \n S(\cdot)B\n_{\ga(L^2(0,t;H),E)} <\infty.$
Since $E$ has finite cotype, $E$ does not contain a copy of $c_0$
and the theorem of Hoffmann-J\o rgensen
and Kwapie\'n implies that $S(\cdot)B \in
\ga(L^2(\R_+;H),E)$.

Thus the implication that remains to be proved is 
(c)$\Rightarrow$(b). A direct proof of (c)$\Rightarrow$(a) would also be of interest, as it would show the equivalence of (a) and (c). 
By standard $H^\infty$-functional methods it is easy to prove that (c) implies the weaker result that $(-A)^{-\alpha}S(\cdot)B$ is in $\ga(L^2(\R_+;H),E)$ for any $\alpha>0$.

Following Weiss  
\cite[Note, page 369]{Weiss:conjectures}, we offer 
$100$ euro for
a positive or negative resolution of these problems.
A consequence of Theorem~\ref{thm:inv} is that the conjecture is true
for bounded and invertible operators $A$ (although this is not of
great practical value) as well as certain
other cases, for instance when $A$ and $B$
diagonalise simultaneously. To see the latter, suppose there is an orthonormal basis $(h_k)_{k\ge 1}$ in $H$ and a sequence $(x_k)_{k\ge 1}$ in $E$
such that $$B h_k = \beta_k x_k, \quad A x_k = -\la_k x_k,$$
with $\la_k>0$ for all $k\ge 1$. Taking $t_k = \la_k$ and assuming the uniform $\ga$--radonification of the set $\big\{\sqrt{\la}R(\la,A)B\suchthat  \la>0\big\}$, we obtain convergence in $E$ of the sum
\begin{align*}
    \EE \Bignorm{ \sum_{k=1}^\infty \ga_k (-A)^{-\einhalb}B h_k}^2 
& = \; \EE \Bignorm{ \sum_{k=1}^\infty \ga_k \la_k^{-\einhalb}\beta_k  x_k}^2 
\\ & =   
    4 \EE \Bignorm{ \sum_{k=1}^\infty \ga_k 
       \frac{t_k^\einhalb  \beta_k}{t_k+\la_k}  x_k}^2 
        =\; 4 \,\EE \Bignorm{ \sum_{k=1}^\infty \ga_k  t_k^\einhalb R(t_k, A)B h_k }^2.
\end{align*}
Consequently, $A^{-\einhalb}B$ is $\ga$--radonifying. 

\medskip
\noindent
{\em Acknowledgment} -- The authors thank the anonymous referee for pointing out a mistake in the original version of the paper.

\def\SUBMITTED{{S}ubmitted}
\def\TOAPPEAR{{T}o appear in }
\def\PREPARATION{{I}n preparation}
\bibliographystyle{ams-pln}
\bibliography{stoch_weiss}

\end{document}